\documentclass[12pt, a4paper, final, oneside, onecolumn, notitlepage]{report}

\usepackage{graphicx}
\usepackage{amsmath,amsfonts}
\usepackage{amssymb}
\usepackage{amsthm}
\usepackage[utf8]{inputenc}
\usepackage{polski}
\usepackage[polish]{babel}
\usepackage[usenames]{color}
\usepackage{enumerate}

\theoremstyle{definition}
\newtheorem{thm}{\small\bf Theorem}

\newcommand{\be}{\begin{eqnarray}}
\newcommand{\ee}{\end{eqnarray}}
\newcommand{\bd}{\begin{displaymath}}
\newcommand{\ed}{\end{displaymath}}
\newcommand{\bee}{\begin{eqnarray*}}
\newcommand{\eee}{\end{eqnarray*}}

\newcommand{\R}{\mathbb{R}}
\newcommand{\N}{\mathbb{N}}

\usepackage{fancyhdr}
\begin{titlepage}
\title{On some geometrical methods leading to martingales useful in measure theory}
\author{Damian Prusinowski}
\date{}
\end{titlepage}

\renewcommand{\maketitle}{

\bibliographystyle{plain}

\newcommand{\linia}{\rule{\linewidth}{0.4mm}}
    \vspace*{1cm}
    \begin{center}\small
    \textsc{Universtity of Lodz}\\
    \textsc{Faculty of Mathematics and Computer Science}\\
		\textsc{ }\\
		\vspace{1.5cm}
    \end{center}
    \vspace{2cm}
    \noindent\linia
    \begin{center}
     \LARGE \textsc{On some geometrical methods leading to martingales useful in measure theory} \\

         \end{center}
     \linia
    \vspace{0.5cm}
    	\begin{center}
					\large \textsc{Damian Prusinowski}
			\end{center}
		\vspace{1cm}
		\vspace{2cm}

    \vspace*{\stretch{6}}
    \begin{center}
    \textsc{Łódź}, 15.09.2017  \textsc{r}.
    \end{center}
}

\begin{document}

\maketitle
\thispagestyle{empty}

\newpage
\thispagestyle{empty}

\newpage

\chapter*{Introduction}
\addcontentsline{toc}{chapter}{Introduction}
In this paper an interval is always an interval being a convex hole $conv\{(0;a),(1;b)\}$ contained in the square $[0;1] \times [0;1] \subset \R^2$. A parallelogram is always a convex hole $P = conv \{(0;a),(0;a+\delta),(1;b),(1;b+\delta)\}$ contained in $[0;1] \times [0;1$. The two-dimensional Lebesgue measure is denoted by $\lambda^2$. The following theorem seems to be unexpected.
\begin{thm}
For some set of intervals $(I_t, t\in T)$ we have \\
$\bigcup_{t \in T} I_t \supset {0} \times [0;1]$ and $\lambda^2(\bigcup_{t \in T} I_t) = 0$.
\end{thm}

This theorem is a rather immediate consequence of the existance of the Nikodym set $N$. Namely for some system $(l_t, t \in T)$ of straight lines on $\R^2$, there exists a system of points $(w_t, t \in T)$ such that $\{w_t, t \in T \} = \R^2$ and $\lambda^2(\bigcup_{t \in T} l_t \setminus \{w_t\}) = 0$. The set $N = \{l_t, t \in T\}$ is known as Nikodym set.\\ \\
\textbf{Proof of Theorem 1.} Let $S \subset T$ be the following subset $S = \{t \in T, w_t \in \{0\} \times [0;1] \}$ and let $w_t = (0;b(t))$ for $t \in S$. We use  affine operations $A_{i,t}(x,y) = (x, b(t) + \frac{1}{i}(y - b(t)))$ for $t \in S$, $i \in \N$ to define the following sets $N_i = \{A_{i,t}[l_t]; t \in S, A_{i,t}[l_t]$ contains an interval$\}$ for $i \in \N$. Let $\{I_{i,t}, t \in S_i\}$ be the set of intervals contained in the lines $A_{i,t}[l_t], t \in S$. It is enough to take $\{ I_t, t \in T' \} = \bigcup_{i \in \N} \{ I_{i,t}, t \in S_i \} \cup \{[0;1] \times \{0\}, [0;1] \times \{1\} \}$. $\square$ \\

Any known construction the Nikodym set is long (see e.g. [1]). It could be interesting that a much shorter proof of Theorem 1. can be obtained by martingale methods. This new construction will contain two naturally separated parts: the geometrical one and the probabilistic one. Our Theorem 2 below is just the geometrical part of this construction. The probabilistic part will be done in [2].

I desire to thank A. Paszkiewicz for pointing me the geometrical problem, discussed in the paper and for valuable discussion.

\newpage
\chapter*{Construction of the Nikodym type set}
\addcontentsline{toc}{chapter}{Construction of the Nikodym set}

Let $\lambda^2$ mean two-dimensional Lebesgue measure. Let parallelogram be a set $P_{b_1,b_2,\delta} 
:= conv\{(0;b_1),(0;b_1+\delta),(1;b_2),(1;b_2+\delta)\}$ for $0 \leq b_1, b_2 \leq 1- \delta$, where $\delta > 0$. Denote also $I_{b,\delta} = \{0\}\times[b, b+\delta]$ for $0 \leq b \leq 1- \delta$, where $\delta > 0$. A interval $I_{b_1,\delta}$ will call the left side of parallelogram $P_{b_1,b_2,\delta}$.

\begin{thm} For all $\epsilon > 0$ there exists parallelograms $P_s$, $s = 1,...,S$ such that
\begin{enumerate}
  \item[(i)] $P_s = P_{b_1^s,b_2^s,\delta^s}$ for $s \in S$ and $\bigcup_{s \in S} I_{b_1^s,\delta^s} \supset \{0\} \times [\epsilon; 1 - \epsilon]$
\end{enumerate}
and for any $x_0, y_0 \in [0,1]$ and $J = [0, x_0]\times[0,1]$ or $[0,1] \times [0,y_0]$ occur the following properties:
\begin{enumerate}
  \item[(ii)] $|\lambda^2(\bigcup_{s \in S} P_s \cap J) - \frac{3}{4} \lambda^2(J)| < \epsilon$;
  \item[(iii)] $|\sum_{s \in S} \lambda^2(P_s \cap J) - \lambda^2(J)| < \epsilon$.
\end{enumerate}
\end{thm}

\textbf{Proof.} Let $n \geq 3$ be an odd number and $Q_i^n = P_{\frac{i}{n^2}; \frac{i+n}{n^2}; n^{-2}}$ i $R_j^n = P_{\frac{j+n}{n^2}; \frac{j}{n^2}; n^{-2}}$, for $i,j \in S(n)$ for $S(n) := \{0,2,...,n^2-n-2\}$. A family of sets $\{P_s, 1 \leq s \leq S \}$ can be obtained in the form $\{Q_i^n, i \in S(n) \}\cup\{R_j^n, j \in S(n) \}$ for sufficiently large $n$.

Property $(i)$ takes place for sufficiently large, odd $n \geq 3$, because the sum of the left sides of the parallelograms \\
$\bigcup_s I_{b_1^s,\delta^s} = \bigcup_{i \in S(n)} \{0\} \times [\frac{i}{n^2}, \frac{i}{n^2} + \frac{1}{n^2}] \cup \bigcup_{j \in S(n)} \{0\} \times [\frac{j + n}{n^2}, \frac{j+n}{n^2} + \frac{1}{n^2}] = \{0\} \times [\frac{1}{n},1 - \frac{1}{n} - \frac{1}{n^2}] \supset \{0\} \times [\epsilon; 1 - \epsilon]$, \\
whenever $n > \frac{2}{\epsilon}$.

We will show that the properties $(ii), (iii)$ occur, only if odd $n$ is greater than certain $n(\epsilon)$, dependent only on $\epsilon$.

\begin{enumerate}
  \item[(ii)] Step $1^o$ Let $J = [0,x_0]\times[0,1]$.
  
$\lambda^2((\bigcup_{i \in S(n)} Q_i^n \cup \bigcup_{j \in S(n)} R_j^n) \cap J) = \lambda^2((\bigcup_{i \in S(n)} Q_i^n \cap J) \cup (\bigcup_{j \in S(n)} R_j^n \cap J)) = \lambda^2(\bigcup_{i \in S(n)} Q_i^n \cap J) + \lambda^2(\bigcup_{j \in S(n)} R_j^n \cap J) - \lambda^2(\bigcup_{i \in S(n)} Q_i^n \cap \bigcup_{j \in S(n)} R_j^n \cap J) = \lambda^2(\bigcup_{i \in S(n)} Q_i^n \cap J) + \lambda^2(\bigcup_{j \in S(n)} R_j^n \cap J) - \lambda^2(\bigcup_{i \in S(n)} \bigcup_{j \in S(n)} Q_i^n \cap R_j^n \cap J) = \Sigma_{i \in S(n)} \lambda^2(Q_i^n \cap J) + \Sigma_{j \in S(n)} \lambda^2(R_j^n \cap J) - \Sigma_{i \in S(n)} \Sigma_{j \in S(n)} \lambda^2(Q_i^n \cap R_j^n \cap J)$, \\
because it is the sum of sets mutually disjoint.

Note that \\  
$\lambda^2(Q_i^n \cap J) = \frac{1}{n^2}x_0$, $\lambda^2(R_j^n \cap J) = \frac{1}{n^2}x_0$ and $\lambda^2(Q_i^n \cap R_j^n) = \frac{1}{2n^3}$ \\
for $i,j = 0, 2,..., n^2-n-2$.

Let $E_1 = \lambda^2(\bigcup_{i \in S(n)} Q_i^n \cap \bigcup_{j \in S(n)} R_j^n \cap [0,\frac{\lfloor n x_0 \rfloor}{n}] \times [\frac{1}{n} - \frac{1}{n^2}, 1 - \frac{2}{n} + \frac{1}{n^2}])$, \\
$E_2 = \lambda^2(\bigcup_{i \in S(n)} Q_i^n \cap \bigcup_{j \in S(n)} R_j^n \cap [0,\frac{\lfloor n x_0 \rfloor}{n}] \times [\frac{1}{n} - \frac{1}{n^2}, 1 - \frac{1}{n}])$ \\
and $E_3 = \lambda^2([\frac{\lfloor n x_0 \rfloor}{n},\frac{\lfloor n x_0 \rfloor + 1}{n}] \times [0, 1] \cup [0,\frac{\lfloor n x_0 \rfloor}{n}] \times [0, \frac{1}{n} - \frac{1}{n^2}] \cup [0,\frac{\lfloor n x_0 \rfloor}{n}] \times [1 - \frac{1}{n}, 1])$.

The following estimate is true \\
$E_1 \leq \Sigma_{i \in S(n)} \Sigma_{j \in S(n)} \lambda^2(Q_i^n \cap R_j^n \cap J) \leq E_2 + E_3$. \\
$\lfloor nx_0 \rfloor \cdot \frac{1}{2n^3} \cdot \frac{n^2 - 3n + 2}{2} \leq \Sigma_{i \in S(n)} \Sigma_{j \in S(n)} \lambda^2(Q_i^n \cap R_j^n \cap J) \leq \lfloor nx_0 \rfloor \cdot \frac{1}{2n^3} \cdot \frac{n^2 - 2n + 1}{2} + \frac{1}{n} + \lfloor nx_0 \rfloor (\frac{1}{n^2} - \frac{1}{n^3}) + \lfloor nx_0 \rfloor \frac{1}{n^2}$. \\
$\frac{\lfloor nx_0 \rfloor}{4} \cdot (\frac{1}{n} - \frac{3}{n^2} + \frac{2}{n^3}) \leq \Sigma_{i \in S(n)} \Sigma_{j \in S(n)} \lambda^2(Q_i^n \cap R_j^n \cap J) \leq \frac{\lfloor nx_0 \rfloor}{4} \cdot (\frac{1}{n} - \frac{3}{2n^2} + \frac{3}{4n^3})$. \\
$\frac{x_0}{4} \cdot (1 - \frac{3}{n} + \frac{1}{n^2}) \leq \Sigma_{i \in S(n)} \Sigma_{j \in S(n)} \lambda^2(Q_i^n \cap R_j^n \cap J) \leq \frac{x_0}{4} \cdot (1 - \frac{3}{2n} + \frac{3}{4n^2})$.

Therefore $|\lambda^2(\bigcup_{i \in S(n)} Q_i^n \cap \bigcup_{j \in S(n)} R_j^n \cap J) - \frac{1}{4}x_0| < \epsilon$, only if $n > \frac{3}{\epsilon}$, and in consequence \\
$|\lambda^2((\bigcup_{i \in S(n)} Q_i^n \cup \bigcup_{j \in S(n)} R_j^n) \cap J) - \frac{3}{4}x_0| < \epsilon$, only if $n > \frac{3}{\epsilon}$.

Step $2^o$ Let $J = [0,1]\times[0,y_0]$, $y_0 \in (0,1)$.

Similarly like in $1^o$ we have \\
$\lambda^2((\bigcup_{i \in S(n)} Q_i^n \cup \bigcup_{j \in S(n)} R_j^n) \cap J) = \Sigma_{i \in S(n)} \lambda^2(Q_i^n \cap J) + \Sigma_{j \in S(n)} \lambda^2(R_j^n \cap J) - \Sigma_{i \in S(n)} \Sigma_{j \in S(n)} \lambda^2(Q_i^n \cap R_j^n \cap J)$

and note that \\
$|\sum_{i \in S(n)} \lambda^2(Q_i^n \cap J) + \sum_{j \in S(n)} \lambda^2(R_j^n \cap J) - y_0| < \epsilon$ only if $n > \frac{2}{\epsilon}$ (see $(iii)$ case $2^o$) and $\lambda^2(Q_i^n \cap R_j^n) = \frac{1}{2n^3}$ for $i,j = 0, 2,..., n^2-n-2$.

Let $E_1 = \lambda^2(\bigcup_{i \in S(n)} Q_i^n \cap \bigcup_{j \in S(n)} R_j^n \cap [0;1] \times [\frac{1}{n} - \frac{1}{n^2}, \frac{\lfloor \frac{1}{2}n^2 y_0 \rfloor}{\frac{n^2}{2}})$, \\
$E_2 = \lambda^2(\bigcup_{i \in S(n)} Q_i^n \cap \bigcup_{j \in S(n)} R_j^n \cap [0;1] \times [\frac{1}{n} - \frac{1}{n^2}, \frac{\lfloor \frac{1}{2}n^2 y_0 \rfloor + 1}{\frac{n^2}{2}}])$ \\
and $E_3 = \lambda^2([0;1] \times [0, \frac{1}{n} - \frac{1}{n^2}])$.

The following estimation is true \\
$E_1 \leq \Sigma_{i \in S(n)} \Sigma_{j \in S(n)} \lambda^2(Q_i^n \cap R_j^n \cap J) \leq E_2 + E_3$; \\
$n \cdot \frac{1}{2n^3} \cdot (\lfloor \frac{1}{2}n^2 y_0 \rfloor - \frac{1}{2}n + \frac{1}{2}) \leq \Sigma_{i \in S(n)} \Sigma_{j \in S(n)} \lambda^2(Q_i^n \cap R_j^n \cap J) \leq n \cdot \frac{1}{2n^3} \cdot (\lfloor \frac{1}{2}n^2 y_0 \rfloor - \frac{1}{2}n + \frac{3}{2}) + \frac{1}{n} - \frac{1}{n^2}$; \\
$\frac{\lfloor \frac{1}{2} n^2 y_0 \rfloor}{2n^2} - \frac{1}{4n} + \frac{1}{4n^2} \leq \Sigma_{i \in S(n)} \Sigma_{j \in S(n)} \lambda^2(Q_i^n \cap R_j^n \cap J) \leq \frac{\lfloor \frac{1}{2} n^2 y_0 \rfloor}{2n^2} - \frac{3}{4n} - \frac{1}{4n^2}$; \\
$\frac{y_0}{4} - \frac{1}{4n} \leq \Sigma_{i \in S(n)} \Sigma_{j \in S(n)} \lambda^2(Q_i^n \cap R_j^n \cap J) \leq \frac{y_0}{4} - \frac{3}{4n} - \frac{1}{4n^2}$.

Therefore $|\Sigma_{i \in S(n)} \Sigma_{j \in S(n)} \lambda^2(Q_i^n \cap R_j^n \cap J) - \frac{1}{4}y_0| < \epsilon$, only if $n > \frac{3}{4\epsilon}$, and in consequence \\
$|\lambda^2((\bigcup_{i \in S(n)} Q_i^n \cup \bigcup_{j \in S(n)} R_j^n) \cap J) - \frac{3}{4}y_0| < \epsilon$, only if $n > \frac{2}{\epsilon}$.

  \item[(iii)] Step $1^o$ Let $J = [0,x_0]\times[0,1]$, $x_0 \in [0,1]$. Note that \\  
$\lambda^2(Q_i^n \cap J) = \frac{1}{n^2}x_0$ and $\lambda^2(R_j^n \cap J) = \frac{1}{n^2}x_0$ for $i,j \in S(n)$.

Therefore

$|\sum_{i \in S(n)} \lambda^2(Q_i^n \cap J) + \sum_{j \in S(n)} \lambda^2(R_j^n \cap J) - \lambda^2(J)| = |\sum_{i \in S(n)} \frac{x_0}{n^2} + \sum_{j \in S(n)} \frac{x_0}{n^2} - x_0| = |\frac{n^2-n}{n^2} x_0 - x_0|= \frac{1}{n}x_0 < \epsilon $ dla $n > \frac{1}{\epsilon}$.

Step $2^o$ Let $J = [0,1]\times[0,y_0]$, $y_0 \in (0,1)$.

Note that \\
$Q_i^n \subset J$ for $i \leq n^2y_0-n$,\\
$Q_i^n \cap J = \emptyset$ for $i \geq n^2y_0$. \\
Similarly \\
$R_j^n \subset J$ for $j \leq n^2y_0-n$,\\
$R_j^n \cap J = \emptyset$ for $j \geq n^2y_0$.

Let $f(n,y_0) = \sum_{i \in S(n)} \lambda^2(Q_i^n \cap J) + \sum_{j \in S(n)} \lambda^2(R_j^n \cap J)$.

$\sum_{i \leq n^2y_0-n , i \in S(n)} \lambda^2(Q_i^n) + \sum_{j \leq n^2y_0-n j \in S(n)} \lambda^2(R_j^n) \leq f(n,y_0) \leq \sum_{i < n^2y_0 , i \in S(n)} \lambda^2(Q_i^n) + \sum_{j < n^2y_0, j \in S(n)} \lambda^2(R_j^n)$.

Because Lebesgue’s measure of each parallelograms is equal $\frac{1}{n^2}$, thus

$\frac{1}{n^2}\frac{\lfloor n^2 y_0 \rfloor - n}{2} \cdot 2 \leq f(n,y_0) \leq \frac{1}{n^2} \frac{\lfloor n^2 y_0 \rfloor - 2}{2} \cdot 2$.

Therefore

$ y_0 - \frac{1}{n} - \frac{1}{n^2} \leq f(n,y_0) \leq y_0 - \frac{2}{n^2}$.

In consequence $|f(n,y_0) - \lambda^2(J)| = |f(n,y_0) - y_0| < \epsilon$, only if $n > \frac{2}{\epsilon}$. $\square$

\end{enumerate}

\newpage
\newpage

\chapter*{References}
\addcontentsline{toc}{section}{References}
[1] M. de Guzman, {\it Differentiation of integrals in $\R^n$}, Springer-Verlag (1975).
\newline
[2] A. Paszkiewicz, D. Prusinowski, {\it On martingale methods and some Nikodym type sets}, in preparation.

\end{document}